\documentclass[10pt,a4paper]{article}
\usepackage{amsmath,amsthm,amssymb,pstricks,pst-node}
\usepackage[dvips]{graphicx}
\parskip=\medskipamount
\usepackage{geometry}

 1   %Caracteres g¢ticos.
 1  %Caracteres "doble palo".
\newcommand{\R}{\mathbb{R}}
\newcommand{\A}{{\cal A}}
\newcommand{\B}{{\cal B}}
\newcommand{\Ro}{{\cal R}}

\newcommand{\e}{\epsilon}
\newcommand{\lag}{\mathfrak{g}}

\def\fpd#1#2{\frac{\partial #1}{\partial #2}}
\newcommand{\F}{\mathbb{F}}

\newcommand{\aff}{{\mathrm{aff}}}
\newtheorem{theorem}{Theorem}
\newtheorem{lemma}{Lemma}

\newtheorem{definition}{Definition}

\theoremstyle{remark}

\title{Routhian reduction for quasi-invariant Lagrangians}
\author{B. Langerock\thanks{e-mail: Bavo.Langerock@architectuur.sintlucas.wenk.be }\\
{\em Sint-Lucas School of Architecture, Hogeschool voor Wetenschap \& Kunst}\\
{\em B-9000 Ghent, Belgium} \\[1.5ex]
F. Cantrijn\thanks{e-mail: Frans.Cantrijn@UGent.be}\mbox{ } and
J. Vankerschaver\thanks{e-mail: Joris.Vankerschaver@UGent.be}
\thanks{Present address: Control and Dynamical Systems, Caltech MC 107--81, Pasadena, CA 911125, USA} \\
{\em Ghent University, Dept. of Mathematical Physics and Astronomy}\\
{\em B-9000 Ghent, Belgium} \ }

\begin{document}
\maketitle
\begin{abstract}
In this paper we describe Routhian reduction as a special case of
standard symplectic reduction, also called Marsden-Weinstein
reduction. We use this correspondence to present a generalization of
Routhian reduction for quasi-invariant Lagrangians, i.e. Lagrangians
that are invariant up to a total time derivative. We show how
functional Routhian reduction can be seen as a particular instance
of reduction of a quasi-invariant Lagrangian, and we exhibit a
Routhian reduction procedure for the special case of Lagrangians
with quasi-cyclic coordinates.  As an application we consider the
dynamics of a charged particle in a magnetic field.
\end{abstract}

\section{Introduction and outline}
In modern geometric approaches to Routhian reduction it is often
mentioned that this reduction technique is the Lagrangian analogue
of symplectic or Marsden-Weinstein reduction~\cite{sympred} (see for
instance the introduction of~\cite{CMR01}). This assertion is
usually justified by the fact that, roughly speaking, for Routhian
reduction one first restricts the system to a fixed level set of the
momentum map and then reduces by taking the quotient with respect to
the symmetry group. In this paper we show, among other things, that
the analogy between Routhian reduction and Marsden-Weinstein
reduction holds at a more fundamental level: in fact we will show
that Routhian reduction is simply a special instance of general
Marsden-Weinstein reduction (from now on referred to as
MW-reduction). More specifically, by applying the MW-reduction
procedure to the tangent bundle of a manifold, equipped with the
symplectic structure induced by the Poincar\'e-Cartan $2$-form
associated with a Lagrangian, we will show that the resulting
reduced symplectic space is `tangent bundle-like', and that the
reduced symplectic structure is again defined by a Poincar\'e-Cartan
form, augmented with a gyroscopic 2-form. Of course, this symplectic
description of the reduced system, obtained via the Routh's
reduction method, is well-known in the literature. The difference
with our approach, however, lies in the fact that we arrive at the
reduced symplectic structure following the Marsden-Weinstein method.
Until now, the symplectic nature of a Routh-reduced system was
obtained either by reducing the variational principle
(see~\cite{jalna00,marsdenrouth} and references therein) or by
directly reducing the second order vector field describing the given
system (see~\cite{mestcram}).

The advantage of interpreting Routhian reduction in terms of
MW-reduction lies in the fact that we are able to extend the concept
of Routhian reduction to quasi-invariant Lagrangian systems, i.e.
Lagrangian systems which are invariant up to a total time
derivative. Such a generalization lies at hand: it is well known
that a quasi-invariant Lagrangian determines a strict invariant
energy and a strict invariant symplectic structure on the tangent
bundle. On the other hand, the actual reduction of quasi-invariant
Lagrangians exploits the full power of MW-reduction and is therefore
in our opinion a very interesting application of this reduction
procedure. The generalization to quasi-invariant Lagrangians is the
main result of this paper.

\paragraph{Lagrangians with a quasi-cyclic coordinate.}
In the remainder of the introduction, we illustrate some of the
concepts used in this paper by means of a simple, but clarifying
example: the case of a Lagrangian with a single quasi-cyclic
coordinate. This is a generalization of the classical procedure of
Routh dealing with Lagrangians with a cyclic coordinate, and will
serve as a conceptual introduction for the geometric techniques
introduced later on, when we deal with the case of general
quasi-invariant Lagrangians in Theorems~\ref{thm:routh}
and~\ref{thm:5}.

We begin by recalling the classical form of Routh's result on the
reduction of Lagrangians with cyclic coordinates (or, stated in a
slightly different way, the reduction of Lagrangians which are
invariant with respect to an abelian group action). For simplicity,
we confine ourselves to the case of one cyclic coordinate.
Subsequently, we will illustrate how this theorem can be extended to
cover the case of quasi-cyclic coordinates.

Given a Lagrangian $L:\R^{2n}\to \R$ for a system with $n$ degrees
of freedom $(q^1,\ldots,q^n)$ for which, say, $q^1$ is a {\em cyclic
coordinate} (i.e. $\partial L/\partial q^1=0$). The momentum
$p_1=\partial L/\partial \dot q^1$ is a first integral of the
Euler-Lagrange equations of motion. If $\partial^2L/\partial\dot
q^1\partial \dot q^1\neq0$ holds, there exists a function $\psi$
such that $p_1=\mu$ is equivalent to $\dot
q^1=\psi(q^2,\ldots,q^n,\dot q^2,\ldots,\dot q^n)$.

\begin{theorem}[Routh reduction~\cite{pars}]
Let $L: \R^{2n} \rightarrow \R$ be a regular Lagrangian for a system with $n$ degrees of freedom $(q^1, \ldots, q^n)$.  Assume that $q^1$ is a cyclic coordinate and that $\partial^2L/\partial\dot
q^1\partial \dot q^1\neq0$ so that $\dot{q}^1$ can be expressed as $\dot{q}^1 = \psi(q^2, \ldots, q^n, \dot{q}^2, \ldots, \dot{q}^n)$. Consider the Routhian $R^\mu:
\R^{2(n-1)}\to \R$ defined as the function $R^\mu = L-\dot q^1\mu$ where
all instances of $\dot q^1$ are replaced by $\psi$. The Routhian is
now interpreted as the Lagrangian for a system with $(n-1)$ degrees
of freedom $(q^2,\ldots,q^n)$.

Any solution $(q^1(t),\ldots,q^n(t))$ of the Euler-Lagrange
equations of motion
\[
\frac{d}{dt}\left(\fpd{L}{\dot q^i}\right)-\fpd{L}{q^i}=0,\
i=1,\ldots,n
\]
with momentum $p_1 = \mu$, projects onto a solution
$(q^2(t),\ldots,q^n(t))$ of the Euler-Lagrange equations
\[
\frac{d}{dt}\left(\fpd{R^\mu}{\dot q^k}\right)-\fpd{R^\mu}{q^k}=0,\
k=2,\ldots,n.
\]
Conversely, any solution of the Euler-Lagrange equations for $R^\mu$ can be lifted to a solution of the Euler-Lagrange equations for $L$ with momentum $p_1 = \mu$.
\end{theorem}
The number of degrees of freedom of the system with Lagrangian
$R^\mu$ is reduced by one, and this technique is called
Routh-reduction. We now formulate a generalization of this theorem
for a Lagrangian system with a \emph{quasi-cyclic coordinate} $q^1$,
i.e. there exists a function $f$ depending on
$(q^1,\ldots,q^n)$ such that
\[
\frac{\partial L}{\partial q^1} = \dot
q^i \frac{\partial f}{\partial q^i}.
\]

If $q^1$ is quasi-cyclic, it is easy  to show that there is an
associated first integral of the Lagrangian system given by $F :=
\partial L/\partial \dot q^1-f$. Note that if
$\partial^2L/\partial\dot q^1\partial \dot q^1\neq 0$, we can again
solve the equation $F = \mu$, where $\mu$ is a constant, to obtain
an expression for $\dot{q}^1$ in terms of the remaining variables.
In the next theorem we now show how the classical procedure of Routh
may be extended to cover the case of a Lagrangian with a
quasi-cyclic coordinate. We defer the proof of this theorem to
section~\ref{ssec:quasicyclic}.

\begin{theorem}[Routh reduction for a quasi-cyclic coordinate]\label{thm:cyclic}
A regular Lagrangian $L:\R^{2n}\to \R$ for a system with $n$ degrees of freedom
$(q^1,\ldots,q^n)$ with a quasi-cyclic coordinate $q^1$ is
Routh-reducible if (i) $\partial^2L/\partial\dot q^1\partial \dot
q^1\neq 0$ and if (ii) there exist $(n-1)$ functions $\Gamma_k$
independent of $q^1$ such that
  \begin{equation} \label{conn}
\fpd{f}{q^k} = \Gamma_k(q^2,\ldots,q^n)\fpd{f}{q^1},\ k=2,\ldots,n .
\end{equation}
For $\mu$ a constant, consider the Routhian $R^\mu:\R^{2(n-1)}\to \R$ defined as
\[
R^\mu= L - (\mu+f(q))(\dot q^1+\Gamma_i\dot q^i),
\]
where all instances of $\dot q^1$ are replaced by the expression
obtained from the equation $\partial L/\partial \dot q^1= \mu +f$.
The Routhian is independent of $q^1$ and can be seen as a Lagrangian
for a system with $(n-1)$ degrees of freedom $(q^2,\ldots,q^n)$.

Then, any solution $(q^1(t),\ldots,q^n(t))$ of the Euler-Lagrange
equations
\[
\frac{d}{dt}\left(\fpd{L}{\dot q^i}\right)-\fpd{L}{q^i}=0,\
i=1,\ldots,n
\]
such that $\partial L/\partial\dot q^1-f = \mu$, projects onto a 
solution $(q^2(t),\ldots,q^n(t))$ of the Euler-Lagrange equations
\[
\frac{d}{dt}\left(\fpd{R^\mu}{\dot q^k}\right)-\fpd{R^\mu}{q^k}=0,\
k=2,\ldots,n.
\]
Conversely, any solution of the Euler-Lagrange equations for $R^\mu$ can be lifted to a solution of the Euler-Lagrange equations for $L$ for which $\partial L/\partial \dot{q}^1 - f = \mu$
\end{theorem}

Readers familiar with methods from differential geometry might
recognize that the functions $\Gamma_k$ determine a
\emph{connection} on the configuration space. The condition (ii)
from the above theorem can be interpreted geometrically as the
existence of a connection for which $df$ annihilates the horizontal
distribution, or alternatively, such that $f$ is covariantly
constant: $Df = 0$ (with $Df$ denoting the restriction of $df$
to the horizontal distribution). It turns out that this condition is
essential to Routhian reduction in the context of quasi-invariant
Lagrangians.

We note that the requirement that $df$ annihilates the horizontal
distribution implies in this case that there exists an equivalent
Lagrangian $L'$ (\emph{i.e.} a Lagrangian that differs from $L$ by a
total time derivative) which is strictly invariant so that Routhian
reduction in the classical sense can be applied. However, we should
warn against dismissing quasi-invariant Routh reduction too hastily
since Routh reduction is possible also for quasi-invariant Lagrangians with nontrivial non-equivariance cocycle. We refer to~\cite{marmo88} for a general discussion on quasi-invariant Lagrangian systems and in particular the property that the vanishing of this non-equivariance cocycle is a necessary condition for a quasi-invariant Lagrangian to be equivalent to a strict invariant Lagrangian.

To conclude this introduction, we note that the study of Routhian
reduction for quasi-invariant Lagrangians was partially inspired on
a technique called \emph{functional Routhian reduction} described
in~\cite{funcrouth}, where it is used to obtain a control law for a
three-dimensional bipedal robot. We will return to this example in
section~\ref{sec:funcrouth}.

\paragraph{Plan of the paper.}
In sections~\ref{sec:clared} and~\ref{sec:routh} we show that
classical Routhian reduction is precisely MW-reduction. We start
with the well-known description of MW-reduction in the cotangent
bundle framework. Although a description of cotangent bundle
reduction may be found in~\cite{cotangentred}, we will elaborate on
this and prove the results because this will show useful when
considering quasi-invariant Lagrangians. Next, in
section~\ref{sec:quasi} we describe MW-reduction for quasi-invariant
Lagrangians. In section~\ref{sec:exam} we conclude with a number of
examples.

\section{Tangent and cotangent bundle reduction}\label{sec:clared}

In this section, we recall some standard results on group actions
and principal bundles and we formulate Marsden-Weinstein reduction
theorem in its standard form. We then specialize to the reduction of
a cotangent bundle with the canonical symplectic form or a tangent
bundle with a symplectic form which is obtained through pullback
along the Legendre transformation.  The material in this section is
well-known and more information can be found in \cite{MarsdenHamRed,
OrRa04}.

\subsection{Momentum maps and symplectic reduction}

\paragraph{Notations.}
Throughout this paper we shall mainly adopt the notations
from~\cite{CMR01} and~\cite{ortegathesis}. Let $M$ be a manifold
on which a group $G$ acts on the right. This action is denoted by
$\Psi: M\times G\to M$ and is such that $\Psi_{gh} = \Psi_h \circ
\Psi_g$ for all $g,h \in G$, with $\Psi_g:\equiv \Psi(\cdot, g)$.
The action $\Psi$ induces a mapping on the Lie-algebra level
\[\varphi: M\times\lag \to TM:(m,\xi)\mapsto \varphi_m(\xi)
=\left.\frac{d}{d\e}\right|_{\e=0}\Psi(m,\exp\e\xi).\] The mapping
$\lag\to \mathfrak{X}(M)$ associating to a Lie-algebra element $\xi$
the corresponding infinitesimal generator
$\xi_M\in\mathfrak{X}(M):m\to \varphi_m(\xi)$ is a Lie-algebra
morphism. The isotropy group $G_m<G$ of an element $m\in M$ is the
subgroup of $G$ determined by $\Psi(m,g)=m$. The Lie-algebra of
$G_m$ is denoted by $\lag_m$. The orbit ${\cal O}_m$ of $m$ is the
subset of $M$ consisting of the elements of the form $\Psi(m,g)$
with $g\in G$ arbitrary. Finally, we will sometimes consider the
dual to $\varphi_m$, i.e. the map $\varphi^*_m: T_m^*M\to \lag^*$.
With a slight abuse of notation, the symbol $\varphi^*$ will also be used to map a 1-form to a $\lag^*$-valued
function on $M$, pointwise defined by
$\varphi^*(\alpha)(m)=\varphi^*_m(\alpha(m))$, with $\alpha$ a 1-form
and $m\in M$ arbitrary.

We will often assume that the action on a manifold $M$ is free and
proper. This guarantees that the space of orbits $M/G$ is a manifold
and that the projection $\pi: M\to M/G$ is a principal fibre
bundle~\cite{koba}. We assume that the reader is familiar with the
concept of associated bundles of a principal manifold and, in
particular, the bundle $\tilde\lag$ associated with the Lie-algebra
$\lag$ on which the group acts on the left by means of the adjoint
action. The adjoint action of $G$ on its Lie-algebra $\lag$ is
denoted by $Ad_g$, and is defined as the differential at the
identity of the conjugation mapping. The dual to the adjoint action
is called the coadjoint action and is denoted by $Ad^*_g$, i.e.
$Ad^*_g(\mu) \in \lag^*$ for $\mu\in\lag^*$. We denote elements in $\tilde\lag$ by $\tilde\xi$ and they represent orbits of points in $Q\times \lag$ under the action of $G$ defined by $(q,\xi) \mapsto (qg,Ad_{g^{-1}}\xi)$ with $q\in Q, g\in G$ and $\xi \in\lag$ arbitrary. In this sense we sometimes write $\tilde\xi = [q,\xi]_G$.

A principal connection on a manifold $M$ on which $G$ acts freely and properly is an
equivariant $\lag$-valued 1-form $\A$ on $M$ such that, in addition,
$\A(\xi_M)=\xi$ for all $\xi\in\lag$. The equivariance property is
expressed by $\A_{\Psi_g(m)}(T\Psi_g(v_m)) = Ad_{g^{-1}}(\A_m(v_m))$, for any
$m\in M, v_m \in T_mM$ and $g\in G$. The kernel of $\A$ determines a
$G$-invariant distribution on $M$ which is called the {\em
horizontal} distribution since it is complementary to the vertical
distribution $V\pi =\ker T\pi$, with $\pi:M\to M/G$. In this paper
we will consider the dual of the linear map $\A_m:T_mM\to \lag$
which is understood to be a map $\A^*_m:\lag^*\to T^*_mM$. If
$\mu\in\lag^*$, then the 1-form $\A^*(\mu):M\to T^*M$ is defined
pointwise by $m\mapsto\A^*_m(\mu)$. Again, with a slight abuse of notation, we sometimes write $\A^*(\mu)=\A_\mu$.

Throughout the paper we encounter products of bundles  over the same base manifold $B$, say $E_1\to B$ and $E_2\to B$. The fibred product $E_1\times_B E_2$ over the base manifold is often denoted simply by $E_1\times E_2$ and consists of pairs $(e_1,e_2)$ with $e_1\in E_1$ and $e_2\in E_2$ such that $e_1$ and $e_2$ project onto the same point in $B$.

\paragraph{Symplectic reduction.}
Let $(M,\omega)$ be a symplectic manifold on which $G$ acts freely on the
right, $\Psi:M\times G\to M$. The action $\Psi$  is {\em canonical}
if $\Psi^*_g\omega=\omega$ for all $g\in G$. If the infinitesimal
generators $\xi_M$ are globally hamiltonian vector fields, i.e. if
there is a function $J_\xi$ for any $\xi\in\lag$ such that
$i_{\xi_M}\omega=-dJ_\xi$, then the map $J:M\to \lag^*$, is called a
momentum map associated to the action.

Following~\cite{Marsden}, we define the {\em non-equivariance
cocycle} associated to a momentum map of the canonical action:
\[
\sigma: G\to \lag^*: g\mapsto J(mg^{-1})-Ad^*_{g^{-1}}(J(m)),
\]
where $m$ is arbitrary in $M$. If $M$ is connected this definition is independent of the
choice of the point $m$ and determines a $\lag^*$-valued one-cocycle $\sigma$
in $G$, i.e.  for $g,h\in G$ it satisfies
\[
\sigma(gh)=\sigma(g)+ Ad^*_{g^{-1}}\sigma(h).
\]
If $M$ is not connected we restrict the analysis to a connected component. Therefore, without further mentioning it, we will always assume that the manifolds we are considering are connected. Given another momentum map $J'$ associated to the same action, its non-equivariance cocycle $\sigma'$ determines the same element as
$\sigma$ in the first $\lag^*$-valued cohomology of $G$, i.e.
$[\sigma]=[\sigma']\in H^1(G,\lag^*)$. Note that for reasons of
conformity, we haven chosen to define $\sigma$
following~\cite{ortegathesis} for left actions: recall that a right
action composed with the group inversion is a left action.

If the moment map is not equivariant one can show (see
\cite{OrRa04}) that it becomes equivariant with respect to the
\emph{affine} action of $G$ on $\lag^*$ determined using the
cocycle $\sigma$ and given by
\[
(g,\mu)\mapsto
Ad^*_{g}\mu +\sigma(g^{-1}).
\]
Due to the fact that $G$ acts freely on $M$ - this is the only case we consider - any value of $J$ is regular and, therefore, $J^{-1}(\mu)$ will be a submanifold of $M$ for all $\mu \in J(M)$~\cite{ortegathesis}.
\begin{theorem}[Marsden-Weinstein reduction]
Let $(M,\omega)$ be a symplectic manifold with $G$ acting freely, properly and
canonically on $M$. Let $J$ be a momentum map for this action with
non-equivariance cocycle $\sigma$. Assume that $\mu \in J(M)$, and denote by $G_\mu$ the isotropy of $\mu$ under the
affine action of $G$ on $\lag^*$. Then $(M_\mu,\omega_\mu)$, with
$M_\mu = J^{-1}(\mu)/G_\mu$, is a symplectic manifold such that the
2-form $\omega_\mu$ is uniquely determined by $i^*_\mu
\omega=\pi_\mu^*\omega_\mu$, with  $i_\mu:J^{-1}(\mu)\to M$ and
$\pi_\mu: J^{-1}(\mu)\to M_\mu=J^{-1}(\mu)/G_\mu$.

Let $H$ denote a function on $M$, which is invariant under the
action of $G$. Then, the Hamiltonian vector field $X_H$ is tangent
to $J^{-1}(\mu)$ and there exists a Hamiltonian $h$ on $M_\mu$ with
$\pi^*_\mu h = i^*_\mu H$, such that the restriction of $X_H$ to
$J^{-1}(\mu)$ is $\pi_{\mu}$-related to $X_h$.
\end{theorem}

\subsection{Cotangent bundle reduction}

Consider now the case of a cotangent bundle $T^\ast Q$ with its
canonical symplectic structure $\omega_Q := d\theta_Q$, where
$\theta_Q$ is the Cartan $1$-form\footnote{Let $\alpha\in T^*Q$,
then $\theta_Q(\alpha)(X)=\langle\alpha,T\pi_Q(X)\rangle$ for
arbitrary $X\in T_\alpha(T^*Q)$.}. Let $G$ be a Lie group
acting freely and properly on $Q$ from the right. Since a cotangent bundle is a special
case of a symplectic manifold, the Marsden-Weinstein theorem
obviously applies to $T^\ast Q$.  However, because of the extra
structure present on a cotangent bundle much more can be said in
this case than one would expect from the Marsden-Weinstein theorem:
see \cite{MarsdenHamRed, cotangentred}.

The group $G$ acts on $Q$ by a right action $\Psi$ and hence also on
$T^\ast Q$ by the cotangent lift of this action: $(g,\alpha)\mapsto
T^*\Psi_{g^{-1}}(\alpha)$. The map $J :=\varphi^*: T^*Q\to \lag^*$,
defined by $\langle J(\alpha_q),\xi\rangle = \langle
\alpha_q,\varphi_q(\xi)\rangle$, is a momentum map for this action.
One can easily show that $J$ is equivariant with respect to the
coadjoint action on $\lag^*$, or in other words, $J \circ
T^*\Psi_{g^{-1}} = Ad_{g}^\ast \circ J$.

Recall that we assume that the action of $G$ is free and
proper so that the quotient $Q/G$ is a manifold.  In this case the
quotient projection $\pi : Q \rightarrow Q/G$ defines a principal
fiber bundle with structure group $G$.  We denote the bundle of
vertical vectors with respect to the projection $\pi$ by $V\pi$. The
subbundle $V^0\pi$ of $T^*Q$ is defined as the annihilator of
$V\pi$.

Fix a principal connection $\A$ on $Q$ and let $\phi^\mu_{\A}$ be the
map $J^{-1}(\mu) \to V^0\pi; \alpha_q \mapsto
\phi^\mu_{\A}(\alpha_q):=\alpha_q-\A^*_q(\mu)$. This is an
equivariant diffeomorphism w.r.t. the standard action of $G_\mu$ on
$V^0\pi$, and its projection onto the quotient spaces is denoted by
$[\phi_\mu^\A]:J^{-1}(\mu)/G_\mu \to V^0
\pi/G_\mu$. The space $V^0\pi$ can be identified with $T^*(Q/G)\times Q$ and, consequently, the quotient space $V^0\pi/G_\mu$ can be identified with the product bundle $T^*(Q/G)\times Q/G_\mu$. We therefore conclude that the choice of a connection $\A$ allows us to identify $J^{-1}(\mu)/G_\mu$ with the bundle $T^*(Q/G)\times Q/G_\mu$ by means of the diffeomorphism $[\phi^\A_\mu]$.

Next, the 1-form $\A_\mu$ (which is also denoted by $\A^*(\mu)$) determines a $G_\mu$-invariant 1-form on $Q$. It is not hard
to show that $d\A_\mu$ is a 2-form on $Q$, projectable to a
2-form $\B_\mu$ on $Q/G_\mu$. This follows from the invariance
under the action of $G_\mu$ and the annihilation of fundamental
vector fields of the form $\xi_Q$ with $\xi$ in the Lie-algebra
$\lag_\mu$ of $G_\mu$. In the following we consider the 2-form on
$T^*(Q/G)\times Q/G_\mu$ determined as the sum of
\begin{itemize}
\item the pull-back to $T^*(Q/G)\times Q/G_\mu$ of $\omega_{Q/G}$ on $T^*(Q/G)$ ;
\item the pull-back to $T^*(Q/G)\times Q/G_\mu$ of $\B_\mu$ on  $Q/G_\mu$.
\end{itemize}
Let $\pi_1, \pi_2$ and $p_\mu$ be the projections
$\pi_1:T^*(Q/G)\times Q/G_\mu \to T^*(Q/G)$, $\pi_2: T^*(Q/G)\times
Q/G_\mu\to Q/G_\mu$ and $p_\mu:Q\to Q/G_\mu$, respectively.  We
further denote the natural injection $V^0\pi\to T^*Q$ by $i_0$. The
above mentioned 2-form on $T^*(Q/G)\times Q/G_\mu$ equals
\[
\pi_1^*\omega_{Q/G} + \pi_2^*\B_\mu.
\]

\begin{figure}[htb]\centering\includegraphics{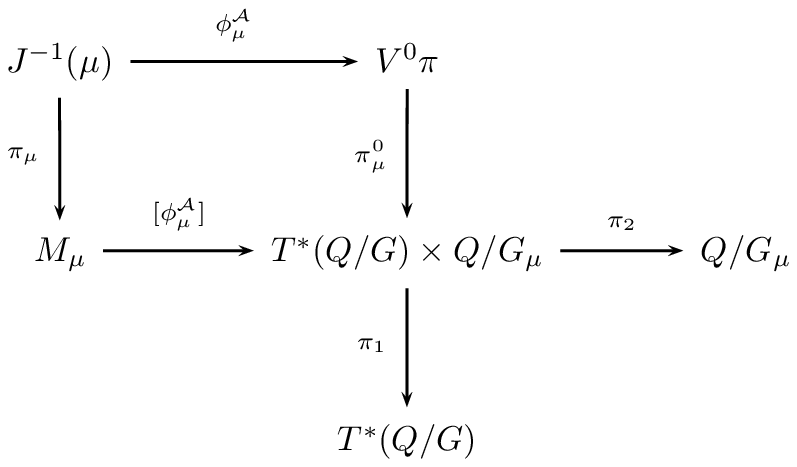}\caption{Cotangent bundle reduction.}\label{fig:diag1}\end{figure}

\begin{theorem}[Cotangent bundle reduction]
Given a free and proper action of $G$ on $Q$ and consider its canonical lift to $T^*Q$. Let $\mu$ be any value of the momentum map, with isotropy subgroup
$G_\mu$. By fixing a principal connection ${\cal A}$, the symplectic
manifold $(M_\mu,\omega_\mu)$ is symplectomorphic to
$(T^*(Q/G)\times Q/G_\mu,\pi_1^*\omega_{Q/G} + \pi_2^*\B_\mu)$, with
symplectomorphism $[\phi^\A_\mu]$.
\end{theorem}
We can summarize this in the diagram presented in Figure~\ref{fig:diag1}.

Although this result is not new and can be found for instance in
\cite{MarsdenHamRed, cotangentred}, we include a proof because its
method will turn out to be useful later on.
\begin{proof} We know that $[\phi_\mu^\A]$ is a diffeomorphism,
and therefore it only remains to show that the symplectic 2-form
$\pi_1^*\omega_{Q/G} + \pi_2^*\B_\mu$ is pulled back to $\omega_\mu$
under this map. We use the fact that $\omega_\mu$ is uniquely
determined by $i^*_\mu\omega_Q=\pi^*_\mu\omega_\mu$, with
$i_\mu:J^{-1}(\mu)\to T^*Q$ the natural inclusion and
$\pi_\mu:J^{-1}(\mu)\to M_\mu$ the projection to the quotient space.
Due to the uniqueness property, it is therefore sufficient to show
that
\begin{equation}\label{eq:symplmorph}\pi^*_\mu([\phi_\mu^\A]^*(\pi_1^*\omega_{Q/G}
+ \pi_2^*\B_\mu))=i^*_\mu\omega_Q.\end{equation}

We will slightly reformulate this condition, by using the fact that
$i^*_\mu\theta_Q = (\phi_\mu^\A)^*\big(i^*_0(\theta_Q
+\pi^*_Q\A_\mu)\big)$ and $[\phi_\mu^\A]\circ\pi_\mu =
\pi^0_\mu\circ\phi_\mu^\A$:
\begin{equation}\label{eq:symplmorph2}
(\pi^0_\mu)^*(\pi_1^*\omega_{Q/G} + \pi_2^*\B_\mu)= d i^*_0(\theta_Q
+\pi^*_Q\A_\mu)
\end{equation}
The latter equality follows easily from the properties of the maps
involved: we have that (i) $(\pi_2\circ\pi^0_\mu)^*\B_\mu =
(\pi_Q\circ i_0)^*d\A_\mu$ and (ii)
$(\pi_1\circ\pi^0_\mu)^*\theta_{Q/G}=i_0^*\theta_{Q}$ hold.
\end{proof}
The above description of cotangent bundle reduction can be seen as
a special case of the more general result stating that if two
symplectomorphic manifolds are both MW-reducible for the same
symmetry group and have compatible actions, then the reduced spaces
are also symplectomorphic. More specifically, given two symplectic
manifolds $(P,\Omega)$ and $(P',\Omega')$ and a symplectomorphism
$f:P\to P'$, i.e. $f^*\Omega'=\Omega$. We assume in addition that
both $P$ and $P'$ are equipped with a canonical free and proper
action of $G$. Let $J:P\to \lag^*$ and $J':P'\to \lag^*$ denote
corresponding momentum maps for these actions on $P$ and $P'$
respectively. We say that $f$ is equivariant if $f(pg)=f(p)g$ for
arbitrary $p\in P$, $g\in G$. Note that the non-equivariance
cocycles for $J$ and $J'$ are equal up to a coboundary. Withouth
loss of generality we assume $f^*J'=J$ and that the non-equivariance
cocycles coincide. This in turn guarantees that the affine actions
on $\lag^*$ coincide and that the isotropy group of an element
$\mu\in\lag^*$ coincides for both affine actions. Finally, fix  a
value $\mu \in \lag^*$ of both $J$ and $J'$.
\begin{theorem}\label{thm:1}
If $f$ is an equivariant symplectic diffeomorphism $P\to P'$ such
that $J'=J \circ f$, then under MW-reduction, the symplectic
manifolds $(P_\mu,\Omega_\mu)$ and $(P'_\mu,\Omega'_\mu)$ are
symplectically diffeomorphic under the map
  \[
  [f_\mu]:P_\mu\to P_\mu';[p]_{G_\mu} \mapsto [f(p)]_{G_\mu}.
  \]
\end{theorem}
\begin{proof}
This is a straightforward result. Since $f$ is a diffeomorphism for
which $J'=J\circ f$, the restriction $f_\mu$ of $f$ to $J^{-1}(\mu)$
determines a diffeomorphism from $J^{-1}(\mu)$ to $J'^{-1}(\mu)$.
The equivariance implies that $f_\mu$ reduces to a diffeomorphism
$[f_\mu]$ from $P_\mu=J^{-1}(\mu)/G_\mu$ to
$P'_\mu=J'^{-1}(\mu)/G_\mu$. It is our purpose to show that
$[f_\mu]^*\Omega'_\mu = \Omega_\mu$ or, since both $\pi_\mu$ and
$\pi'_\mu$ are projections, that
$\pi_\mu^*\Omega_\mu=f^*_\mu(\pi'^*_\mu\Omega'_\mu)$. The
determining property for $\Omega_\mu$ and $\Omega'_\mu$ is
$\pi^*_\mu\Omega_\mu = i^*_\mu\Omega$ (similarly for $\Omega'_\mu$).
From diagram chasing we have that $i^*_\mu\Omega = f_\mu^*(i'^*_\mu
\Omega')$. Then \[\pi_\mu^*\Omega_\mu =
i^*_\mu\Omega=f^*_\mu(i'^*_\mu \Omega') = f^*_\mu(\pi'^*_\mu
\Omega_\mu') =\pi^*_\mu([f_\mu]^*\Omega_\mu'),\] since
$\pi'_\mu\circ f_\mu = [f_\mu]\circ\pi_\mu$ by definition. This
concludes the proof.
\end{proof}

\subsection{Tangent bundle reduction}
We start by recalling the symplectic formulation of Lagrangian
systems on the tangent bundle $TQ$ of a manifold $Q$, and its
relation to the canonical symplectic structure on $T^*Q$ through the
Legendre transform. Next, we shall consider Lagrangians invariant
under the action of $G$, and study a general Marsden-Weinstein
reduction scheme for such systems.

\begin{definition} A Lagrangian system is a pair $(Q,L)$ where $Q$ is
called the configuration manifold and $L$ is a smooth function on $TQ$.
A Lagrangian system $(Q,L)$ is said to be regular if the fibre derivative
$\F L:TQ\to T^*Q; v_q\mapsto \F L(v_q)$ is a diffeomorphism.
The map $\F L$ is called the Legendre transformation and is defined by
\[
\langle \F L(v_q),w_q\rangle =
\left.\frac{d}{d\epsilon}\right|_{\epsilon=0}L(v_q+\epsilon w_q),
\]
for arbitrary $v_q,w_q\in T_qQ$.
\end{definition}
\begin{definition} Given a free and proper action $\Psi$ of $G$ on $Q$,
then a Lagrangian system $(Q,L)$ is said to be invariant if $L$ is an
invariant function for the lifted action $(v_q,g)\mapsto T\Psi_g(v_q)$.
\end{definition}

\begin{figure}[h]\centering
\includegraphics{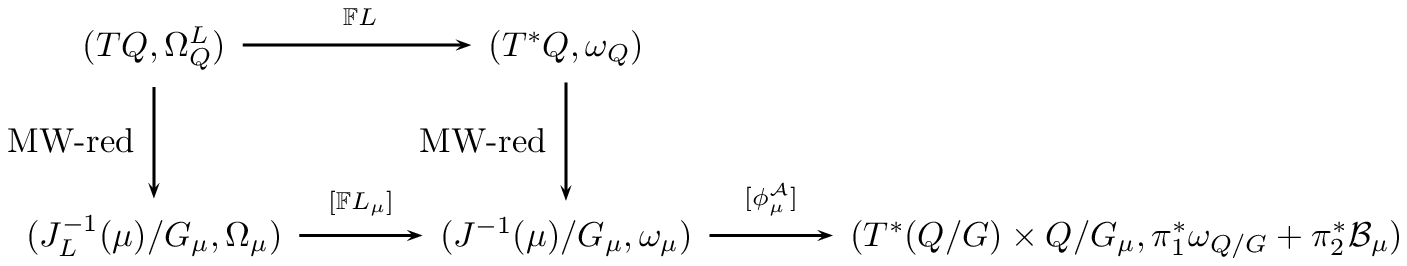}\caption{Diagram relating tangent and cotangent reduction}\label{fig:diag2}\end{figure}

Given a regular Lagrangian system $(Q,L)$, one can define a
symplectic structure on $TQ$ by using the Legendre transform: we
denote the 2-form on $TQ$ obtained by pulling back $\omega_Q$
under $\F L$, by $\Omega^L_Q=(\F L)^*\omega_Q$. We will only
consider regular Lagrangians throughout this paper. The following
results are standard.
\begin{theorem} The lifted action $T\Psi$ of $G$ on $TQ$ is a canonical
action for the symplectic manifold $(TQ,\Omega^L_Q)$. A momentum map is
given by $J_L=J\circ\F L:TQ\to \lag^*$, and $J_L$ is equivariant w.r.t. to
the coadjoint action on $\lag^*$. Furthermore the Legendre transformation
is an equivariant symplectomorphism between the symplectic manifolds
$(TQ,\Omega^L_Q)$ and $(T^*Q,\omega_Q)$.
\end{theorem}
The above theorem guarantees that Theorem~\ref{thm:1} is applicable.
We are now ready to draw the diagram in Figure~\ref{fig:diag2}, with
$\mu\in\lag^*$.

Next, we will show that the manifold $J^{-1}_L(\mu)/G_\mu$ is
diffeomorphic to the fibred product $T(Q/G)\times Q/G_\mu$ if $L$
satisfies an additional regularity assumption. Lagrangians
satisfying this condition are called {\em $G$-regular}. We shall
compute the map $[\phi^\A_\mu]\circ [\F L_\mu]$ and show that it
coincides with a Legendre transform for a function defined on
$T(Q/G)\times Q/G_\mu$. This fact will eventually allow us to show
that the reduced symplectic spaces are again originating from a
Lagrangian system on $Q/G_\mu$. The Lagrangian of this `reduced'
Lagrangian system is precisely the Routhian known from classical
Routhian reduction.

We use the fixed connection $\A$ on $Q$ to identify $TQ/G$ with
the bundle $T(Q/G)\times \tilde\lag$ in the standard way. This identification is obtained as follows: let $[v_q]_G\in TQ/G$ be arbitrary and fix a representative $v_q\in TQ$. The image in $T(Q/G)\times\tilde\lag$ of $[v_q]_G$ is defined as the element $(T\pi(v_q),\tilde\xi)$ with $\pi:Q\to Q/G$ and $\tilde\xi = [q,\A(v_q)]_G \in\tilde\lag$. This map is invertible and determines a diffeomorphism (see for instance~\cite{CMR01}). To define the inverse: let $(v_x,\tilde\xi)$ be arbitrary in $T(Q/G)\times\tilde\lag$, and consider the tangent vector  $v_q = (v_x)^h_q + \varphi_q(\xi)$ at $q\in \pi^{-1}(x)$, with $(v_x)^h_q$ the horizontal lift determined by $\A$ and $\xi$ such that $\tilde\xi=[q,\xi]_G$. The inverse of $(v_x,\tilde\xi)$ is the orbit $[v_q]_G\in TQ/G$ (the latter is well defined: one can show that it is independent of the point $q$, see also~\cite{CMR01}).

Completely analogous one can show that $TQ/G_\mu$ is diffeomorphic to $T(Q/G)\times Q/G_\mu\times \tilde\lag$. Indeed, let $[v_q]_{G_\mu} \in TQ/G_\mu$ be arbitrary and fix a representative $v_q\in TQ$, then the image of $[v_q]_{G_\mu}$ is defined by $(v_x,p_\mu(q),\tilde\xi)$, with $T\pi(v_q)=v_x$ and $\tilde\xi=[q,\A(v_q)]_G$ (recall that $p_\mu:Q\to Q/G_\mu$).  The construction of the inverse map uses the previous diffeomorphism and consists of three steps. Let $(v_x,y,\tilde\xi) \in T(Q/G)\times Q/G_\mu\times \tilde\lag$ be arbitrary. First, we consider the element $[v_q]_G$ in $TQ/G$ which is the inverse of $(v_x,\tilde\xi)\in T(Q/G)\times\tilde\lag$. Secondly, we take a representative $v_q$ of $[v_q]_G$ at a point $q\in p_\mu^{-1}(y)$. And finally, we consider $[v_q]_{G_\mu}$. It is not hard to show that this inverse is well-defined (i.e. independent of the chosen representative $v_q$).

An invariant Lagrangian $L$ determines a function on the quotient $TQ/G$, and
under the identification determined above, a function $l$ on
$T(Q/G)\times\tilde\lag$. We define the fibre derivative $\F_{\tilde\xi} l:T(Q/G)\times \tilde\lag\to T(Q/G)\times\tilde\lag^*$ by
\[
\langle \F_{\tilde\xi} l(v_x,\tilde\xi),(v_x,\tilde\eta)\rangle: =
\left.\frac{d}{d\epsilon}\right|_{\epsilon = 0}l(v_x,\tilde\xi
+\epsilon\tilde\eta)\,.
\]
\begin{definition} An invariant Lagrangian $L$ is said to be
$G$-regular if the map
$\F_{\tilde\xi} l:T(Q/G)\times\tilde\lag\to T(Q/G)\times\tilde\lag^*$
is a diffeomorphism.
\end{definition}

We remark here that according to the previous definition, $G$-regularity depends on the chosen connection $\A$. However, we mention here that $G$-regularity can alternatively be defined as a condition on $L$ directly. We refer the reader to~\cite{BM} for a detailed discussion on $G$-regularity.

A momentum value $\mu $ determines in the quotient spaces a mapping $\tilde\mu:Q/G_\mu \to \tilde\lag^*$ as follows: let $y\in Q/G_\mu$ be arbitrary
\[
\langle\tilde\mu(y),\tilde\xi\rangle = \langle \mu,\xi\rangle,
\]
with $\xi$ the unique representative of $\tilde \xi=[q,\xi]_G$ at a
point $q\in p_\mu^{-1}(y)$. Recall that $p_\mu$ denotes the projection $p_\mu:Q\to Q/G_\mu$. Due to the identification $TQ/G_\mu\cong T(Q/G)\times Q/G_\mu\times \tilde\lag$ the manifold $J^{-1}_L(\mu)/G_\mu$ is a subset of $T(Q/G)\times Q/G_\mu\times \tilde\lag$. In the following lemma we characterize this subset in terms of $\tilde\mu$ and $\F_{\tilde\xi} l$.

\begin{lemma}
There is a one-to-one correspondence between $J^{-1}_L(\mu)/G_\mu$
and the subset of $T(Q/G)\times Q/G_\mu\times\tilde\lag$ determined
as the set of points $(v_x,y,\tilde\xi)$ that satisfy the condition $\F_{\tilde\xi}l(v_x,\tilde\xi)=(v_x,\tilde\mu(y))$.
\end{lemma}
\begin{proof}

Consider a point $[v_q]_{G_\mu}$ in $J^{-1}_L(\mu)/G_\mu$ and let
$v_q$ be a representative. Then, by definition, $J_L(v_q)=\mu$, and
since $L$ is invariant, we have $L(v_q)=l(v_x,\tilde\xi)$, with
$(v_x,\tilde\xi)$ the element in $TQ/G\cong T(Q/G)\times
\tilde\lag$ corresponding to $[v_q]_G$.
Using the definition of the momentum map $J_L$, we obtain
\[
\langle J_L(v_q),\eta\rangle =\left. \frac{d}{d\epsilon}
\right|_0L(v_q+\epsilon \varphi_q(\eta)) = \left. \frac{d}{d\epsilon}
\right|_0 l(v_x, \tilde\xi+\epsilon \tilde\eta) = \langle
\F_{\tilde\xi}l(v_x,\tilde\xi),(v_x,\tilde\eta)\rangle,
\]
with $\tilde\eta=[q,\eta]_G$.
\end{proof}

\begin{lemma}\label{lem:identi}
Let $L$ be $G$-regular invariant Lagrangian. Then there is a
diffeomorphism between $J_L^{-1}(\mu)/G_\mu$ and $T(Q/G)\times
Q/G_\mu$.
\end{lemma}
\begin{proof}
We define a map $J^{-1}_L(\mu)/G_\mu\to T(Q/G)\times Q/G_\mu$ and
its inverse. Let $[v_q]_{G_\mu}\in J^{-1}_L(\mu)/G_\mu$, with
$v_q\in J^{-1}_L(\mu)$ a representative at $q$. We again use the
fixed connection $\A$ on $Q$, and we introduce the maps:
  \begin{align*}
  & p_1([v_q]_{G_\mu}):= T\pi(v_q), &&
   p_2([v_q]_{G_\mu}) :=p_\mu(q), &&
   p_3([v_q]_{G_\mu}):=[q,\A(v_q)]_G,
  \end{align*}
with $\pi:Q\to Q/G$ and $p_\mu: Q\to Q/G_\mu$. These maps $p_{1,2,3}$ are simply the restrictions to $J^{-1}_L(\mu)/G_\mu$ of the projections on the first, second and third factors in the product $T(Q/G)\times Q/G_\mu \times \tilde\lag$. It is easily verified that $(p_1,p_2): J^{-1}_L(\mu)/G_\mu\to T(Q/G)\times Q/G_\mu$ is smooth.

We now define the inverse map $\psi_\mu$ of $(p_1,p_2)$. Let
$(v_x,y)\in T(Q/G)\times Q/G_\mu$ be arbitrary and define the element $\tilde
\xi\in\tilde\lag_x$ such that $(v_x,\tilde\xi) = (\F_{\tilde\xi}l)^{-1}(v_x,
\tilde\mu(y))$ (here we use the condition that $L$ is $G$-regular).
Now consider the tangent vector $v_q = (v_x)^h_q + \varphi_q(\xi)$,
where $\xi$ is such that $\tilde\xi=[q,\xi]_G$ and $q\in
p_{\mu}^{-1}(y)$. By construction we have on the one hand that
$J_L(v_q)=\mu$ and on the other hand
$(p_1,p_2)([v_q]_{G_\mu})=(v_x,y)$.
\end{proof}
Combined with Figure~\ref{fig:diag2}, we can now draw the diagram in Figure~\ref{fig:diag3} below.
\begin{figure}[htb]\centering
\includegraphics{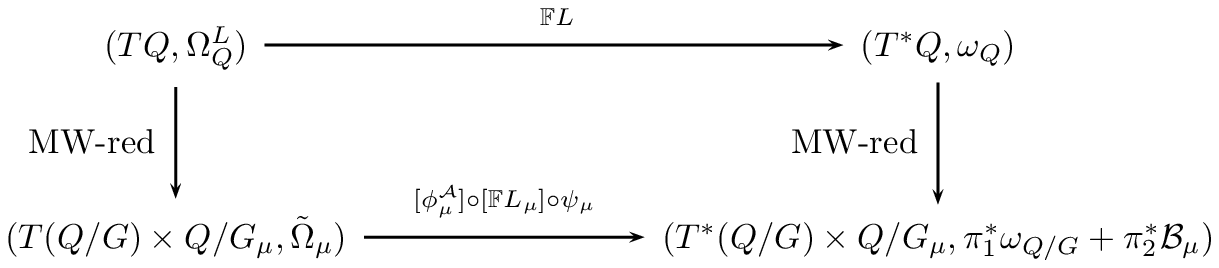}\caption{Diagram relating tangent and cotangent reduction for $G$-regular Lagrangians}\label{fig:diag3}\end{figure}

There is an interesting local criterium for a Lagrangian $L$ to be
$G$-regular. We say that $L$ is \emph{locally $G$-regular} if given
a point $[v_q]_{G_\mu}$ in $J^{-1}_L(\mu)/G_\mu$, there is a
neighborhood $U$ of $[v_q]_{G_\mu}$ such that the restriction $(p_1,p_2)|_U$ is a
diffeomorphism from $U$ to its image $(p_1,p_2)(U)$.

\begin{lemma}\label{lem:trivtang}
An invariant Lagrangian $L$ is locally $G$-regular if one of the
following two equivalent conditions hold:
\begin{enumerate}
\item\label{decomp}
$T(J^{-1}_L(\mu))\oplus V_{J^{-1}_L(\mu)}\varphi=T_{J^{-1}_L(\mu)}(TQ)$,
with $V\varphi \subset V\tau_Q \subset T(TQ)$ defined as the set of
tangent vectors of the form $(\varphi_q(\xi))^v_{v_q}=\xi_Q^v(v_q)$,
$\xi \in \lag$ arbitrary, where $\xi_Q$ is the fundamental vector field of the action on $Q$ corresponding to $\xi$, and $\cdot ^v$ denotes the vertical lift $TQ\times TQ\to T(TQ)$.
\item The `vertical' Hessian of $l$, defined as
\[
D^2l(v_x,
\tilde\xi)(\tilde\eta,\tilde\eta'):=\fpd{^2l}{\epsilon\partial\epsilon'} (v_x,\tilde\xi+\epsilon\tilde\eta+\epsilon'\tilde\eta')|_{\epsilon=\epsilon'=0}
\]
for $v_x\in T(Q/G)$ and
$\tilde\xi,\tilde\eta,\tilde\eta'\in\tilde\lag_x$, is invertible.
\end{enumerate}
\end{lemma}
\begin{proof}
Note that for all $v_q \in J^{-1}_L(\mu)$,
\[
    \dim T_{v_q}(J^{-1}_L(\mu)) + \dim V\varphi(v_q) =
    (\dim TQ - \dim \lag) + \dim \lag =
    \dim T_{v_q}(TQ).
\]
The direct-sum decomposition in~(\ref{decomp}) is therefore
equivalent to the statement that $T(J^{-1}_L(\mu))\cap V_{J^{-1}_L(\mu)}\varphi =0$. We will now prove that this is equivalent to the vertical Hessian of $l$ being invertible.

Assume that the intersection $T(J^{-1}_L(\mu))\cap V\varphi$
contains a non-zero element.  Such an element is necessarily of the
form $(\xi_Q)^v(v_q)$, where $\xi \in \lag$ and $\xi \ne 0$.
Expressing the fact that this element is contained in
$T(J^{-1}_L(\mu))$ implies that for every $\eta \in \lag$, $\langle
TJ_L((\xi_Q)^v(v_q)), \eta \rangle  = 0$.  This can be made more
explicit as follows:
\begin{align*}
\langle TJ_L((\xi_Q)^v(v_q)), \eta \rangle & =
\frac{d}{ds} J_L(v_q + s \varphi_q(\xi))(\eta) \Big|_{s = 0} \\
& = \frac{d}{ds} \F_{\tilde \xi} l(v_x, \tilde{\zeta}
+ s \tilde{\xi}) (\tilde{\eta}) \Big|_{s = 0} \\
& = D^2 l(v_x, \tilde{\zeta})(\tilde{\xi},
\tilde{\eta}),
\end{align*}
where we have decomposed $v_q$ in its vertical and horizontal parts
as $v_q =(v_x)^h_q + \varphi_q(\zeta)$. Since this holds for every
$\eta \in \lag$, we conclude that $(\xi_Q)^v(v_q)$ is contained in
the intersection $T(J^{-1}_L(\mu))\cap V_{J^{-1}_L(\mu)}\varphi$ if and only if the
associated section $\tilde{\xi}$ is in the null space of $D^2 l(v_x, \tilde{\zeta})$.  Hence, the two statements in
Lemma~\ref{lem:trivtang} are equivalent.

If $D^2 l(v_x, \tilde{\zeta})$ is invertible, then via the
implicit function theorem, the reduced Legendre transformation is
locally invertible. The method of proof of the previous
Lemma~\ref{lem:identi} can be used to show that, locally,
$(p_1,p_2)$ is invertible.
\end{proof}
Note in passing that if the given $L$ is a mechanical Lagrangian, i.e. it is of type kinetic minus potential, then $L$ is $G$-regular if the locked inertia tensor, defined by the restriction of the kinetic energy metric to the fundamental vector fields (see e.g. \cite{marsdenrouth}), is
non-degenerate. In our language, the `reduced' locked inertia
tensor coincides with $D^2 l_{x}: \tilde\lag\to
\tilde\lag^*$ in the following sense:
\[
D^2l_{x} (\tilde\xi,\tilde\eta) = J(q)(\xi,\eta),
\]
with $\tilde\xi=[q,\xi]_G$ and $\tilde\eta=[q,\eta]_G$ arbitrary.

\section{Routhian reduction}\label{sec:routh}

In this section, we make a start with Routhian reduction. We
consider a Lagrangian $L : TQ \rightarrow \mathbb{R}$ which is
invariant under the action of a Lie group $G$ and as before we consider a connection
$\mathcal{A}$ in the bundle $\pi: Q \rightarrow Q/G$.  Furthermore,
let $\mu \in \mathfrak{g}^\ast$ be a fixed momentum value, and
define the function $R^\mu$ as $R^\mu=L-\A_\mu$ (recall that
$\A_\mu:TQ\to \R$ is the connection 1-form contracted with
$\mu\in\lag^*$).  By definition $R^\mu$ is $G_\mu$-invariant and
in particular, its restriction to $J^{-1}_L(\mu)$ is reducible to a
function $[R^\mu]$ on the quotient $J^{-1}_L(\mu)/G_\mu$. In turn,
we denote the function on $T(Q/G)\times Q/G_\mu$ corresponding to
$[R^\mu]$ by ${\cal R}^\mu$, i.e.${\cal R}^\mu=\psi_\mu^*[R^\mu]$.
The function ${\cal R}^\mu$ is called the {\em Routhian}.

We begin by reconsidering some aspects from the reduction theory of
tangent bundles, which we relate to the geometry of the Routhian.
Recall from the diagram in Figure~\ref{fig:diag3} that we may write
the symplectic 2-form $\tilde\Omega_\mu$ obtained from
MW-reduction as
\[
\left([\phi^\A_\mu]\circ[\F
L_\mu]\circ\psi_\mu\right)^*\left(\pi_1^*\omega_{Q/G} +
\pi_2^*\B_\mu\right).
\]

\begin{lemma}\label{lem:legendre}
The map $[\phi^\A_\mu]\circ[\F L_\mu]\circ\psi_\mu$ is the fibre
derivative of the Routhian ${\cal R}^\mu$, i.e. for
$(v_x,y),(w_x,y)\in T(Q/G)\times Q/G_\mu$ arbitrary
\[
\left\langle \big([\phi^\A_\mu]\circ[\F
L_\mu]\circ\psi_\mu\big)(v_x,y), (w_x,y)\right\rangle =
\left.\frac{d}{d\e}\right|_{\e=0}{\cal R}^\mu(v_x+\e
w_x,y)=:\langle\F {\cal R}^\mu(v_x,y),(w_x,y)\rangle.
\]
\end{lemma}
\begin{proof}
Fix elements $(v_x,y)\in T(Q/G)\times Q/G_\mu$ and fix a $v_q\in
J^{-1}_L(\mu)$ that projects onto $\psi_\mu(v_x,y)$. By definition
of the maps involved, we have
\[
\big([\phi^\A_\mu]\circ[\F L_\mu]\circ\psi_\mu\big)(v_x,y) =
(\pi^0_\mu \circ \phi^\A_\mu)(\F L(v_q)) = \pi^0_\mu(\F
L(v_q)-\A_\mu(q)).
\]
Fix a curve $\e\mapsto \zeta(\e)$ in $J^{-1}_L(\mu)$ that projects
onto the curve $\e\mapsto\psi_\mu(v_x+\e w_x,y)$ in
$J^{-1}_L(\mu)/G_\mu$ and such that $\zeta(0)=v_q$ and $\dot
\zeta(0)$ is vertical to the projection $\tau_Q\circ
i_\mu:J^{-1}_L(\mu)\to Q$. The existence of such a curve is best
shown using Lemma~\ref{lem:trivtang} and some coordinate
computations. For that purpose, fix a bundle adapted coordinate
chart on $Q\to Q/G$, and let $(x^i,g^a)$ denote the coordinate
functions with $i=1,\ldots,\dim Q/G$ and $a=1,\ldots,\dim G$. From
Lemma~\ref{lem:trivtang}, where it was shown that $TJ^{-1}_L(\mu)$ is transversal to $V\varphi$, we deduce that $(x^i, v^i, g^a)$ are (local)
coordinate functions for $J^{-1}_L(\mu)$, with $(x^i,v^i)$ a
standard coordinate chart on $T(Q/G)$ associated to $(x^i)$ on
$Q/G$. In this coordinate chart we put $v_q=(x^i_0,v^i_0,g^a_0)$ and $w_x=(x^i_0,w^i_0)$,
and we define the curve $\zeta(\epsilon)$ to be the curve
$\epsilon\mapsto (x^i_0,v^i_0+\epsilon w^i_0,g^a_0)$. Then the
tangent to $\zeta$ at $\e=0$ is the vertical lift of some $w_q\in
T_qQ$ with $T\pi(w_q)=w_x$.

Finally, from the definition of ${\cal R}^\mu$,
\[
\left.\frac{d}{d\e}\right|_{\e=0}{\cal R}^\mu(v_x+\e w_x,y) =
\left.\frac{d}{d\e}\right|_{\e=0}\left(L(\zeta(\e)) -
\A_\mu(\zeta(\e))\right) = \langle \F L(v_q)-\A_\mu(q),w_q\rangle.
\]
Since $\F L(v_q)-\A_\mu(q) \in V^0\pi$, the right-hand side of this
equation can be rewritten as a contraction with $(w_x,y)$:
\[
\left.\frac{d}{d\e}\right|_{\e=0}{\cal R}^\mu(v_x+\e w_x,y)  =
\langle \pi^0_\mu(\F L(v_q)-\A_\mu(q)), (w_x,y)\rangle.
\]
This concludes the proof.
\end{proof}
The above lemma allows us to compute the reduced symplectic 2-form
on the manifold $T(Q/G)\times Q/G_\mu$.
\[
(\F {\cal R}^\mu)^*\big( \pi_1^*\omega_{Q/G} + \pi_2^*\B_\mu\big) =
(\F {\cal R}^\mu)^*\big( \pi_1^*\omega_{Q/G}\big) +
\overline\pi_2^*\B_\mu,
\]
with $\overline\pi_2:T(Q/G)\times Q/G_\mu \to Q/G_\mu$. In order to
complete the symplectic reduction we now study the energy-function
(this is the Hamiltonian function for the Euler-Lagrange equations).
Recall that the energy $E_L$ corresponding with the Lagrangian
system $(Q,L)$ is the function on $TQ$ defined by $E_L(v_q)=\langle
\F L(v_q),v_q\rangle - L(v_q)$, for $v_q\in TQ$ arbitrary. The
energy for the Routhian ${\cal R}^\mu$ is defined by
\[
E_{{\cal R}^\mu} (v_x,y) = \langle\F {\cal
R}^\mu(v_x,y),(v_x,y)\rangle- {\cal R}^\mu(v_x,y),
\]
with $(v_x,y)\in T(Q/G)\times Q/G_\mu$ arbitrary.
\begin{lemma}\label{lem:energy}
The energy $E_{{\cal R}^\mu}$ is the reduced Hamiltonian, i.e. it
satisfies:
\[
((p_1,p_2)\circ\pi_\mu)^*E_{{\cal R}^\mu} = i^*_\mu E_L,
\]
with $\pi_\mu: J^{-1}_L(\mu)\to J^{-1}_L(\mu)/G_\mu$ and
$i_\mu:J^{-1}_L(\mu)\to TQ$.
\end{lemma}
\begin{proof}
Let $v_q\in J^{-1}_L(\mu)$, such that
$((p_1,p_2)\circ\pi_\mu)(v_q)=(v_x,y)$. Then
\begin{align*}
i^*_\mu E_L(v_q)&= \langle \F L(v_q),v_q\rangle-L(v_q)\\
&= \langle  \phi^\A_\mu(\F L_\mu(v_q)) + \A^*_q(\mu),v_q\rangle - L(v_q)\\
&= \left\langle \big([\phi^\A_\mu]\circ[\F
L_\mu]\circ\psi_\mu\big)(v_x,y),(v_x,y)\right\rangle - {\cal R}^\mu (v_x,y).
\end{align*}
Using Lemma~\ref{lem:legendre} this concludes the
proof.
\end{proof}

We end this section with some additional definitions in order to
interpret the MW-reduced system as a Lagrangian system (we also
refer to~\cite{BM}). For that purpose consider a manifold $M$ fibred
over $N$ with projection $\kappa:M\to N$. Roughly said, a Lagrangian
$L$ with configuration space $M$ is said to be {\em intrinsically
constrained} if it does not depend on the velocities of the fibre
coordinates of $\kappa:M\to N$.  This is made more precise in the
following definition.

\begin{definition}
A Lagrangian system $(M,L)$ on a fibred manifold $\kappa:M\to N$ is
intrinsically constrained if $L$ is the pull-back of a function $L'$
on $T_MN=TN\times_N M$ along the projection $TM\to T_MN$.
\end{definition}
For notational simplicity we will identify $L$ with $L'$. If we fix
a coordinate neighborhood $(x^i ,y^a)$ on $M$ adapted to the fibration, we
can write the Euler-Lagrange equations for this system. The fact
that the Lagrangian is intrinsically constrained is locally expressed by the
fact that $L(x,\dot x,y)$ is independent of $\dot y$, and the
Euler-Lagrange equations then read:
\[
\frac{d}{dt}\left(\fpd{L}{\dot x^i}\right)-\fpd{L}{x^i}=0,\
i=1,\ldots,n \mbox{ and }\fpd{L}{y^a}=0,\ a=1,\ldots,k.
\]
The latter $k$ equations determine constraints on the system.
We now wish to write these equations as Hamiltonian equations
w.r.t. a presymplectic 2-form on $T_MN$. For that purpose, we associate
to the Lagrangian $L:T_MN\to \R$  a {\em Legendre transform} $\F L :T_MN\to T_M^*N$.
The definition is given by, for $(v_n,m),(w_n,m)\in T_MN$ arbitrary
\[
\langle \F L (v_n,m),(w_n,m)\rangle = \left.\frac{d}{d\epsilon}\right|_{\epsilon=0}
L(v_n+\epsilon w_n,m),
\]
In coordinates $\F L(x^i ,\dot x^i,y^a)$ simply reads $(x^i,\partial
L/\partial \dot x^i,y^a)$. Finally, if we write the projection $T^*_MN\to
T^*N; (\alpha_n,m)\to \alpha_n$ by $\kappa_1$, then it is not hard
to show that the pull-back to $T_MN$ of the canonical
symplectic form $\omega_{N}$ under the map $\kappa_1\circ\F
L:T_MN\to T^*N$ determines a presymplectic 2-form, locally equal to
\[
d\left(\fpd{L}{\dot x^i}\right)\wedge dx^i.
\]
We define the energy as the function
\[
E_L:T_MN\to \R; (v_n,m)\to \langle \F L(v_n,m),(v_n,m)\rangle-L,
\]
and the solutions $m(t)$ to the Euler-Lagrange equations solve the equation
\[
\left.\left(i_{\dot\gamma}(\kappa_1\circ\F L)^*\omega_N = -dE_L\right)\right|_{\gamma}
\]
with $\gamma(t)=(\dot n(t),m(t))$ and $n(t)=\kappa(m(t))$ (see also~\cite{gotaya,gotayb}).

If the original intrinsically constrained Lagrangian system $(M,L)$
is non-conservative with a {\em gyroscopic force term}, i.e. a
2-form $\beta$ on $M$ is given and the force term is the function $TM\to T^*M; v_m\mapsto -i_{v_m}\beta_m$, then the Euler-Lagrange equations of motion are Hamiltonian w.r.t
(pre)-symplectic form $(\kappa_1\circ\F L)^*\omega_N
+\kappa_2^*\beta$ and with Hamiltonian $E_L$:
 \[
\left.\left(i_{\dot \gamma}\left((\kappa_1\circ\F L)^*\omega_N +\kappa_2^*\beta\right)= -dE_L\right)\right|_{\gamma}.
\] Here $\kappa_2$
denotes the projection to the second factor in $T_MN$, i.e.
$\kappa_2:T_MN\to M$. In the case of Routhian reduction, the reduced space is of this type: the total space corresponds to $Q/G_\mu$ and the base space $N$ to $Q/G$.
\begin{theorem}\label{thm:routh}
Given a $G$-invariant, $G$-regular Lagrangian $L$ defined on the
configuration space $Q$. Then the MW-reduction of the symplectic manifold
$(Q,\Omega_L)$ for a momentum value $J_L=\mu$ is the symplectic manifold
\[ \big(T(Q/G)\times Q/G_\mu, (\F {\cal R}^\mu)^*\big(
\pi_1^*\omega_{Q/G}\big) + \overline\pi_2^*\B_\mu\big).\] The
reduced Hamiltonian of $E_L$ is the energy $E_{\Ro^\mu}$. The
equations of motion for this Hamiltonian vector field are precisely
the Euler-Lagrange equations of motion for an intrinsically
constrained Lagrangian system on $Q/G_\mu\to Q/G$ with Lagrangian
$\Ro^\mu$ and gyroscopic force term determined by the 2-form
$\B_\mu$ on $Q/G_\mu$.
\end{theorem}
It is remarkable that the 2-form $\B_\mu$ is such that the presymplectic 2-form $(\F {\cal R}^\mu)^*\big( \pi_1^*\omega_{Q/G}\big) + \overline\pi_2^*\B_\mu\big)$ is symplectic. A next step in Routhian reduction would be to identify $\B_\mu$ as a
2-from which is built up o.a. out of the curvature of $\A$ and a nondegenerate part on the fibres of $Q/G_\mu\to Q/G$. Since this is
not the scope of this paper, we refer the reader
to~\cite{BM,marsdenrouth}.

\section{Quasi-invariant Lagrangians}\label{sec:quasi}
In this section we study a possible generalization of the Routhian
reduction procedure to quasi-invariant Lagrangians. We refer the
reader to~\cite{marmo88} and references therein for further details
on quasi-invariant Lagrangians. We assume throughout this section
that $Q$ is a connected manifold, which ensures that given a
function $f$ for which $df=0$ implies that $f$ is constant.

\subsection{Quasi-invariance and cocycles}

We begin by defining what it means for a Lagrangian to be quasi-invariant
under a group action.  We then show that the transformation behaviour of
a quasi-invariant Lagrangian induces a certain cocycle on the space of
1-forms, and we study the properties of this cocycle.

\begin{definition}
A Lagrangian system $(Q,L)$ is quasi-invariant if the Lagrangian
satisfies
\[
(T\Psi_g)^* L (v_q)  = L(v_q) + \langle v_q , dF_g(q)\rangle ,
\]
with $v_q$ arbitrary and for some function $F: G\times Q\to \R$. We denote a quasi-invariant Lagrangian system as a triple $(Q,L,F)$.
\end{definition}
Clearly, the function $F$ is not arbitrary: from the
fact that $\Psi$ defines a right action it follows that
$(T\Psi_{gh})^*L=\big((T\Psi_g)^*\circ (T\Psi_h)^*\big) L$ and one can see that $dF: G\to {\cal X}^*(Q)$ should define a group
$1$-cocycle with values in the $G$-module of 1-forms on $Q$, i.e.
for $g_1,g_2\in G$ arbitrary
\[
\Psi^*_{g_1}dF_{g_2} -dF_{g_1g_2} + dF_{g_1} =0.
\]
Consider the
map $f: \lag \times Q\to \R$ defined by \[f(\xi,q) =\left.
\frac{d}{d\e}\right|_{\e=0} F(\exp\e\xi,q).\] Clearly, $f$ is linear
in its first argument, and thus determines a map $Q\to \lag^*$ which
is denoted by the same symbol. We now define a $1$-cocycle with values in $\lag^*$.

\begin{lemma} The map
 \[
\sigma_F:G\to \lag^* : g\mapsto Ad^*_{g^{-1}} f(q)-
\big(\Psi^*_{g^{-1}}f\big)(q) +
Ad_{g^{-1}}^*\big(\varphi_q^*\big(dF_{g^{-1}}(q)\big)\big).
\]
does not depend on the chosen point $q$ and determines a group
1-cocycle with values in $\lag^*$.
\end{lemma}
\begin{proof}
We first show that the differential of \[q\mapsto f_{Ad_g\xi}(q)-
\Psi^*_gf_\xi(q) + \langle (Ad_g\xi)_Q(q),dF_g(q)\rangle\]vanishes
for arbitrary $\xi\in\lag$. This implies that the above definition
of $\sigma_F$ does not depend on the chosen point $q$.

We start from the cocycle property of the map $g\mapsto dF_g$, i.e.
we have $\Psi^*_{g_1}dF_{g_2} -dF_{g_1g_2} + dF_{g_1} =0$. Let
$g_1=g$ and $g_2=\exp\e\xi$, and take the derivative at $\e=0$, then
\[
\Psi^*_g df_\xi - \left.\frac{d}{d\e}\right|_{\e=0}dF_{g\exp \e
\xi}=0.
\]
To compute the second term we again use the cocycle property with
$g_1=g(\exp\e\xi) g^{-1},g_2=g$, i.e. $dF_{g\exp\e\xi} = dF_{(\exp\e
Ad_{g}\xi )g}=\Psi^*_{\exp\e Ad_{g}\xi}dF_g +dF_{\exp\e Ad_{g}\xi}$.
The derivative with respect to $\epsilon$ at $0$ equals
\[
\left.\frac{d}{d\e}\right|_{\e=0}dF_{g\exp \e \xi}(q) =
d\left(\langle (Ad_g\xi)_Q,dF_{g}\rangle\right)(q) + df_{Ad_g\xi}.\]
We conclude that the map $\langle\sigma_F(g),\xi\rangle =
f_{Ad_{g^{-1}}\xi}(q)- \Psi^*_{g^{-1}}f_\xi(q) + \langle
(Ad_{g^{-1}}\xi)_Q(q),dF_{g^{-1}}(q)\rangle$ is independent of $q$
and therefore well-defined. From straightforward computations it
follows that it is a group 1-cocycle with values in $\lag^*$: for
$g_1,g_2$ arbitrary,
  \[  Ad^*_{g_1^{-1}}\sigma_F(g_2) - \sigma_F(g_1g_2)+ \sigma_F(g_1)=0.\]
  This concludes the proof.
\end{proof}
This  $1$-cocycle induces a $\lag^*$-valued
1-cocycle on the Lie-algebra, given by
\[
\xi \mapsto -ad^*_\xi f + \xi_Q(f) -\varphi^*(df_\xi);
\]
and hence also
a real valued $2$-cocycle $\Sigma_f(\xi,\eta)=
\xi_Q(f_\eta)-\eta_Q(f_\xi)-f_{[\xi,\eta]}$.
This is the cocycle used in the infinitesimal version of quasi-invariant
Lagrangians discussed in for instance~\cite{marmo88}. If only an
infinitesimal action is given, i.e. a Lie algebra morphism $\lag \to \mathfrak{X}(Q);
\xi\mapsto \xi_Q$; or by complete lifting, an infinitesimal action on
$TQ$, then the above definition of 1-cocycle $\sigma_F$ corresponds
infinitesimally to $\Sigma_f$. It is often easier to compute
$\Sigma_f$ instead of $\sigma_F$ in examples, see
section~\ref{sec:exam}.

\subsection{The momentum map}
As mentioned in the introduction, Noether's theorem is applicable to
quasi-invariant Lagrangians as well: for each Lagrangian that is
quasi-invariant under a group action, there exists a momentum map which is conserved.  In this section, we study the properties of
this momentum map, with a view towards performing symplectic
reduction later on.

We begin by investigating the equivariance of the Legendre
transformation.

\begin{lemma}\label{lem:quasilegendre}
Let $(Q,L,F)$ denote a quasi-invariant system. Then, for $g\in G$ arbitrary, the Legendre map
$\F L$ transforms as
\[
\F L(T\Psi_g(v_q))= T^*\Psi_{g^{-1}}\left(\F L(v_q) +
dF_{g}(q)\right)= T^*\Psi_{g^{-1}}\left(\F L(v_q)\right) -
dF_{g^{-1}}(qg).
\]
\end{lemma}
\begin{proof}To show this equality, fix an element $w_{qg}\in TQ$,
and let $w_q = T\Psi_{g^{-1}}(w_{qg})$. Then, by definition of the
fibre derivative,
\begin{align*}
\langle w_{qg} , \F L(T\Psi_g(v_q))\rangle &=
\left.\frac{d}{d\e}\right|_{\e=0} L(T\Psi_g(v_q)+\e w_{qg})\\&
=\left.\frac{d}{d\e}\right|_{\e=0} \left( L(v_q+\e w_q) + \langle
v_q +\e w_q,dF_g(q)\rangle\right)\\&= \langle w_{qg},
T^*\Psi_{g^{-1}}\left(\F L(v_q) + dF_{g}(q)\right)\rangle.
\end{align*}
From $\Psi^*_g(dF_{g^{-1}}) = -dF_g$ (let $g_1=g, g_2=g^{-1}$ in the
cocycle identity for $dF$) we have the property that
$T^*\Psi_{g^{-1}}(dF_g(q))=- dF_{g^{-1}}(qg)$ for $q\in Q$ and $g\in
G$ arbitrary. This concludes the proof.
\end{proof}

The above lemma justifies the next definition.
\begin{definition}
Let $(Q,L,F)$ denote a quasi-invariant Lagrangian system. Then we
define a right action $\Psi_\aff$ on $T^\ast Q$ as follows. For
$\alpha_q\in T^*Q$ arbitrary, we put:
\[
\Psi_{\aff, g}(\alpha_q)= T^*\Psi_{g^{-1}}\left(\alpha_q +
dF_{g}(q)\right) =T^*\Psi_{g^-1}\left(\alpha_q\right) -
dF_{g^{-1}}(qg).
\]
We say that $\Psi_\aff$ is the affine action on $T^*Q$ associated to
the 1-cocycle $dF$.
\end{definition}
We should check that the affine action is well defined. For that
purpose, we need to verify that for $g_1,g_2$ arbitrary
\[
T^*\Psi_{(g_1g_2)^{-1}}\big(\alpha_q+dF_{g_1g_2}(q)\big) =
T^*\Psi_{g_2^{-1}}\left(T^*\Psi_{g_1^{-1}}\big(\alpha_q+dF_{g_1}(q)\big)
+ dF_{g_2}(qg_1)\right).
\]
This is a straightforward consequence from the fact that $dF$ is a
group 1-cocycle.
\begin{lemma}
Let $(Q,L,F)$ denote a quasi-invariant Lagrangian system. Then,
\begin{enumerate}
\item\label{lem:quasitota} the lifted action $T\Psi$ is a canonical action
for the symplectic structure $(TQ,\Omega^L_Q)$;
\item\label{lem:quasitotb} the map $J^f_L=\varphi^*\circ\F L-\tau^*_Qf: TQ\to \lag^*$
is a momentum map with non-equivariance cocycle $\sigma_F$ and the energy
$E_L$ is an invariant function on $TQ$;
\item\label{lem:quasitotc} the affine action $\Psi_\aff$ is a canonical
action for the symplectic structure $(T^*Q,\omega_{Q})$; the map
$J^f=\varphi^*-\pi^*_Q f$ is a momentum map with non-equivariance cocycle $\sigma_F$;
\item\label{lem:quasitotd} $\F L$ is a symplectomorphism between
$(TQ,\Omega^L_Q)$ and $(T^*Q, \omega_Q)$, and is equivariant w.r.t.
to the lifted action on $TQ$ and the affine action on $T^*Q$ associated to $dF$.
\end{enumerate}
\end{lemma}
\begin{proof}
The affine action $\Psi_\aff$ on $T^*Q$ acts by symplectic transformations, i.e. from local computations it follows that
\[
(\Psi_{\aff, g})^*\theta_Q = \theta_Q + \pi_Q^* dF_g.
\]Together with Lemma~\ref{lem:quasilegendre}, i.e. $\F L \circ T\Psi_g = \Psi_{\aff,g} \circ \F L$, assertions (\ref{lem:quasitota}) and (\ref{lem:quasitotd}) follow:
\[
(T\Psi_g)^*\Omega^L_Q = d (\F L \circ T\Psi_g)^* \theta_Q =  \F
L^*  d\Psi_{\aff,g}^*\theta_Q = \Omega_L.
\]
The latter equality holds since $\theta_Q$ is invariant under the
affine action up to an exact form.

To show that $J^f_L$ is a momentum map we use an argument involving
coordinate expressions. Let $(q^i), i=1,\ldots,\dim Q$ denote
coordinate functions on $Q$, and let $(q^i,\dot q^i)$ be the associated
coordinate system on $TQ$. Then it is not hard to show that
\[
\xi_{TQ}\left(\fpd{L}{\dot q^i}\right) =
\fpd{f_\xi}{q^i}-\fpd{\xi^j_Q}{q^i}\fpd{L}{\dot q^j},
\]
holds, with $j=1,\ldots,\dim Q$ and $\xi_Q^j$ the coordinate
expression of $\xi_Q$: $\xi_Q= \xi_Q^j\partial_j$. From some tedious
computations it follows that
\[
i_{\xi_{TQ}}\Omega^L_Q = -dJ^f_\xi
\]
for $\xi\in\lag$ arbitrary.  We now compute the non-equivariance
cocycle of $J^f_L$. Fix any $\xi \in\lag$ and $v_q\in T_qQ$, then
\begin{align*}
\langle J^f_L(T\Psi_g(v_q)),\xi\rangle &
= \langle \F L(T\Psi_g(v_q)),\varphi_{qg}(\xi) \rangle - f_\xi(qg)\\
&=\left\langle T^*\Psi_{g^{-1}}\big(\F L(v_q) +dF_g(q)\big),
T\Psi_g\big(\varphi_q(Ad_g\xi)\big)\right\rangle -f_\xi(qg)\\
&= \langle \F L(v_q),\varphi_q(Ad_g\xi)\rangle -f_{Ad_g\xi}(q)+
\big(f_{Ad_g\xi}(q) -f_\xi(qg)+\langle dF_{g}(q) ,
\varphi_q(Ad_g\xi)\rangle\big)\\&= \langle Ad^*_g
J^f_L(v_q),\xi\rangle + \langle\sigma_F(g^{-1}),\xi\rangle.
\end{align*}
Finally, the fact that the energy is invariant easily follows from
Lemma~\ref{lem:quasilegendre}, and from this we conclude
that~(\ref{lem:quasitotb}) holds.

Since $\F L$ is a
symplectic diffeomorphism and since $J^f\circ\F L= J^f_L$, we
conclude that $J^f$ is a momentum map with cocycle $\sigma_F$. This
proves~(\ref{lem:quasitotc}).
\end{proof}
The above lemma ensures that the equivariance conditions for
Theorem~\ref{thm:1} are satisfied. In that case we can study the
MW-reduction and the structure of the corresponding quotient spaces.
If these quotient spaces are `tangent and cotangent bundle like' we
shall say that the MW-reduction is a Routhian reduction procedure.

Following Theorem~\ref{thm:1} we have that the reduced Legendre
transformation $[\F L_\mu]$ is a symplectic diffeomorphism relating
the symplectic structures on $(J^f_L)^{-1}(\mu)/G_\mu$ and
$(J^f)^{-1}(\mu)/G_\mu$. The subgroup $G_\mu$ is the isotropy
subgroup of the affine action of $G$ on $\lag^*$ corresponding to
the 1-cocycle $\sigma_F$. We now study the structure of the reduced
manifolds $(J^f_L)^{-1}(\mu)/G_\mu$ and $(J^f)^{-1}(\mu)/G_\mu$, and
their respective symplectic 2-forms.

Let $\A$ be a principal connection with horizontal projection
operator $TQ\to TQ : v_q\mapsto v_q^h :=v_q-\varphi_q(\A_q(v_q))$.
Similarly, we can restrict a covector $\alpha_q$ to horizontal
tangent vectors: $T^*Q\mapsto T^*Q : \alpha_q\mapsto \alpha_q^h$,
with $\langle v_q,\alpha_q^h\rangle = \langle
v_q^h,\alpha_q\rangle$. Note that $\alpha_q^h = \alpha_q
-(\A^*_q\circ\varphi^*_q)(\alpha_q)$. The covariant exterior derivative~(see \cite{koba}) of a function $\lambda$ on $Q$ is denoted by $D\lambda$ and is defined pointwise as $D \lambda_q = d\lambda_q^h$. We first study the symplectic
structure of $(J^f)^{-1}(\mu)/G_\mu$. Similar to the invariant
situation, we contract the connection 1-form on the Lie-algebra
level with $\mu+f$ to obtain a 1-form $\A_\mu^f = q\mapsto
\langle\mu+f(q), \A_q\rangle$ on $Q$.

\begin{lemma} Consider a quasi-invariant Lagrangian system $(Q,L,F)$,
for which there exists a principal connection $\A$ such that $DF_g=0$,
for arbitrary $g\in G$. Then,
\begin{enumerate}
\item\label{lem:quasireda} the 2-form $d\A_\mu^f$ is invariant under
the action of $G_\mu$ on $Q$ and is projectable to a 2-form on $Q/G_\mu$
denoted by $\B_\mu^f$;
\item\label{lem:quasiredb} there exists a symplectic diffeomorphism
\[[\phi_\mu^{\A,f}]: ((J^f)^{-1}(\mu)/G_\mu,\omega_\mu) \to
(T^*(Q/G)\times Q/G_\mu,\pi_1^*\omega_{Q/G} + \pi_2^*\B^f_{\mu}).\]
\end{enumerate}
\end{lemma}
\begin{proof} For the proof of both statements we rely on the following
identities, for $g\in G_\mu$ and $q\in Q$:
\begin{align*}
& \A^*_{qg}=T^*\Psi_{g^{-1}}\circ \A^*_q\circ Ad_{g^{-1}}^*\\
&\mu=Ad^*_{g^{-1}}\mu+\sigma_F(g)\\
& Ad^*_{g^{-1}} f(qg) = f(q)- (Ad^*_{g^{-1}}\circ\varphi^*_{qg} )(dF_{g^{-1}}(qg)) + \sigma_F(g)\\
& Ad^*_{g^{-1}}\circ\varphi^*_{qg} = \varphi^*_q\circ T^*\Psi_{g}.
\end{align*}
\ref{lem:quasireda}. The first statement is proven if we can show
that $\A_\mu^f$ is invariant under $G_\mu$ up to an exact 1-form.
Thus consider any element $q\in Q$ and $g\in G_\mu$, then
\begin{align*}
(\Psi^*_g\A^f_\mu)(q)&= \langle \mu+f(qg),Ad_{g^{-1}}\cdot \A_q\rangle \\
&= \left\langle(\mu - \sigma_F(g)) + \big(f(q)-
(Ad^*_{g^{-1}}\circ\varphi^*_{qg} )(dF_{g^{-1}}(qg)) + \sigma_F(g)\big), \A_q  \right\rangle\\
&= (\A^f_\mu)(q) + dF_g(q).
\end{align*}
The latter equality  holds because $dF_g^h(q)=DF_g(q)=0$.
To show that the 2-form is projectable, we prove in addition that
$i_{\xi_Q}d\A^f_\mu=0$. This follows on the one hand from ${\cal
L}_{\xi_Q}\A^f_\mu = df_\xi$ which is obtained using the previous
equation with $g=\exp \e\xi$, and on the other hand from ${\cal
L}_{\xi_Q} = i_{\xi_Q}d+di_{\xi_Q}$:
\[
i_{\xi_Q}d\A^f_\mu = {\cal L}_{\xi_Q}\A^f_\mu - d f_\xi = 0.
\]
\ref{lem:quasiredb}. Similar to the case of an invariant Lagrangian
system we relate $(J^f)^{-1}(\mu)$ with $V^0\pi$ by means of the
connection: $\phi_\mu^{\A,f}: (J^f)^{-1}(\mu)\to
V^0\pi;\alpha_q\mapsto \alpha_q - \A^*_q(\mu+f(q))$. The next step
is to study the affine action of $G_\mu$ on $(J^f)^{-1}(\mu)$
through this diffeomorphism. Let $g\in G_\mu$, and $\alpha_q\in
(J^f)^{-1}(\mu)$, then
\begin{align*}
\phi_\mu^{\A,f}\big( T^*\Psi_{g^{-1}}(\alpha_q+dF_g(q))\big) &=
T^*\Psi_{g^{-1}}(\alpha_q+dF_g(q)) - \A^*_{qg}(\mu+f(qg))\\&=
T^*\Psi_{g^{-1}}\big(\alpha_q - \A^*_{q}(\mu+f(q)\big)+ dF_g^h(qg),
\end{align*}
We conclude that $\phi_\mu^{\A,f}$ is equivariant w.r.t the affine
action on $(J^f)^{-1}(\mu)$ and the standard lifted action on $T^*Q$
restricted to $V^0\pi$ if the condition $DF_g=0$ holds. The
reduced map is denoted by $[\phi_\mu^{\A,f}]$ and maps
$(J^f)^{-1}(\mu)/G_\mu$ to $T^*(Q/G)\times Q/G_\mu$.  The fact that
it is a symplectic map the follows from analogous arguments as in
the invariant case.
\end{proof}

\subsection{The reduced phase space}

We are now ready to take the final step towards a Routhian reduction
procedure for quasi-invariant Lagrangians. It concerns the
realization of $(J^f_L)^{-1}(\mu)/G_\mu$ as a tangent space
$T(Q/G)\times Q/G_\mu$. We therefore reintroduce $G$-regular
quasi-invariant Lagrangians. It should be clear that the definitions
here are also valid in the strict invariant case. Let $R^\mu =
L-\A^f_\mu$ denote the `Routhian' as a function on $TQ$. We first
show that it is $G_\mu$-invariant. For that purpose let $g\in G_\mu$
and $v_q\in T_qQ$, then
\begin{align*}
R^\mu(T\Psi_g(v_q)) &= L(v_q) + \langle dF_g(q),v_q\rangle - \langle
 (\Psi^*_g\A^f_\mu)(q),v_q\rangle \\ & = L(v_q) -\langle
\A^f_\mu(q),v_q\rangle =R^\mu(v_q).
\end{align*}
We know from the strict invariant case that $TQ/G_\mu$ can be
identified with $T(Q/G)\times Q/G_\mu\times \tilde\lag$. Let us
denote $\Re^\mu$ denote the function on the latter space obtained
from projecting $R^\mu$. We now define the fibre derivative
$\F_{\tilde\xi} \Re^\mu$ of $\Re^\mu$ w.r.t the $\tilde\lag$-fibre:
\[
\langle \F_{\tilde\xi} \Re^\mu(v_x,y,\tilde\xi),(v_x,y,\tilde\eta)\rangle =
\left.\frac{d}{d\epsilon}\right|_{\epsilon=0} \Re^\mu(v_x,y,\tilde\xi+\epsilon\tilde\eta),
  \]
with $(v_x,y,\tilde\xi)\in T(Q/G)\times Q/G_\mu\times
\tilde\lag$ and $\tilde\eta\in\tilde\lag_x$ arbitrary.
\begin{definition}
Let $(Q,L,F)$ denote a quasi-regular Lagrangian system. We say that
the system is $G$-regular if the function $\F_{\tilde\xi}
\Re^\mu:T(Q/G)\times Q/G_\mu\times \tilde\lag\to T(Q/G)\times
Q/G_\mu\times \tilde\lag^*$ is a diffeomorphism.
\end{definition}
It is not so hard to show that there is a one-to-one identification
with $(J^f_L)^{-1}(\mu)/G_\mu$ and the set of points $(v_x,y,\tilde\xi)$ in $T(Q/G)\times Q/G_\mu\times \tilde\lag$ for which
$\F_{\tilde\xi} \Re^\mu (v_x,y,\tilde\xi)= (v_x,y,0)$. We consider the map $(p_1,p_2):(J^f_L)^{-1}(\mu)/G_\mu\to T(Q/G)\times
Q/G_\mu$ taking a point $[v_q]_G$ to the first two factors of the corresponding point $(T\pi(v_q),p_\mu(q),[q,\A(v_q)]_G)$ in the fibred product $T(Q/G)\times Q/G_\mu\times \tilde\lag$.

\begin{lemma}
If $(Q,L,F)$ is a $G$-regular quasi-invariant Lagrangian system,
then the mapping $(p_1,p_2):(J^f_L)^{-1}(\mu)/G_\mu\to T(Q/G)\times
Q/G_\mu$ is a diffeomorphism with inverse $\psi_\mu$.
\end{lemma}
The proof is completely analogous to the proof of
Lemma~\ref{lem:identi}: the inverse of $(v_x,y)$ is defined as the point in $(J^f_L)^{-1}(\mu)/G_\mu$ that corresponds to $(\F_{\tilde\xi}\Re^\mu)^{-1}(v_x,y,0)$ in $T(Q/G)\times Q/G_\mu\times\tilde\lag$. Let $[R^\mu]$ denote the quotient of the
restriction of $R^\mu$ to $(J^f_L)^{-1}(\mu)$. Similar to the
previous case we define $\Ro^\mu$ to be function on $T(Q/G)\times
Q/G_\mu$ such that $(p_1,p_2)^*(\Ro^\mu)=[R^\mu]$. Note that
$\Ro^\mu$ could also be obtained by $\Ro^\mu(v_x,y)=
\Re^\mu(v_x,y,\tilde\xi)$, with $(v_x,y,\tilde\xi)=
(\F_{\tilde\xi}\Re^\mu)^{-1}(v_x,y,0)$.

\begin{lemma}\label{thm:4}
Let $(Q,L,F)$ denote a $G$-regular quasi-invariant Lagrangian system
and let $\A$ be a principal connection such that $DF_g=0$. Let
$\mu$ denote a value of $J^f_L$. Then
\begin{enumerate}
\item\label{thm:4a} the map $[\phi_\mu^{\A,f}]\circ [\F L_\mu]\circ \psi_\mu$
is the fibre derivative of $\Ro^\mu$;
\item\label{thm:4b} the energy of $\Ro^\mu$ is the MW-reduced hamiltonian
of the energy $E_L$ on the symplectic manifold $(Q,\Omega^L_Q)$.
\end{enumerate}
\end{lemma}
The proof is again completely similar to the proof of
Lemma's~\ref{lem:legendre} and~\ref{lem:energy}. We conclude that
the MW-reduction of a $G$-regular quasi-invariant Lagrangian $L$ is
again a `Lagrangian' system on the manifold $T(Q/G)\times Q/G_\mu$,
with Lagrangian $\Ro^\mu$: the symplectic structure is of the form
$(\F \Ro^\mu)^*(\pi_1^*\omega_{Q/G}) + \overline\pi_2^*\B_\mu^f$.
\begin{theorem}\label{thm:5}
Let $(Q,L,F)$ denote a $G$-regular quasi-invariant Lagrangian system
and let $\A$ be a principal connection such that $DF_g=0$. Let
$\mu$ denote a value of $J^f_L$.  Then the MW-reduction of
the symplectic manifold $(Q,\Omega_L)$ for the regular momentum
value $\mu$ is the symplectic manifold \[ \big(T(Q/G)\times Q/G_\mu,
(\F {\cal R}^\mu)^*\big( \pi_1^*\omega_{Q/G}\big) +
\overline\pi_2^*\B^f_\mu\big).\] The reduced Hamiltonian is the
energy $E_{\Ro^\mu}$. The equations of motion for this Hamiltonian
vector field are precisely the Euler-Lagrange equations of motion
for an intrinsically constrained Lagrangian system on $Q/G_\mu\to
Q/G$ with Lagrangian $\Ro^\mu$ and gyroscopic force term associated
to the 2-form $\B^f_\mu$ on $Q/G_\mu$.
\end{theorem}

\section{Examples}\label{sec:exam}

\subsection{Quasi-cyclic coordinates}\label{ssec:quasicyclic}

We continue here the description started in the introduction of a
Lagrangian $L$ with a single quasi-cyclic coordinate. Recall that if
$(q^1,\ldots,q^n)$ are coordinates on $Q = \R^n$ and $L(q^i,
\dot{q}^i)$ is a Lagrangian, then we say that $q^1$ is quasi-cyclic
if there exists a function $f(q^1, \ldots, q^n)$ such that
\[
\fpd{L}{q^1}=\dot q^i\fpd{f}{q^i}.
\]
The group $G = \R$ acts on $\R^n$ by translation in $q^1$. Since
$\lag \equiv \R$, a principal connection $\A$ here becomes an
ordinary $G$-invariant $1$-form on $\R^n$. The infinitesimal version
$\Sigma_f$ of the definition of the cocycle $\sigma_F$ is
identically zero, and we can conclude that also $\sigma_F$ vanishes.
Since the group is abelian, we have that $G_\mu =G$. The quotient
space is $T(Q/G)$ and $Q/G$ is labeled by the configuration space
coordinates $(q^2,\ldots,q^n)$.

The condition that the system should be $G$-regular is locally
expressed by $\partial^2L/\partial \dot q^1\partial \dot q^1\neq0$
and, secondly, the condition that there exists a (principal)
connection $\A$ such that $Df=0$ (i.e. $df$ restricted to the
horizontal distribution should vanish) boils down to the condition
that there should exist functions $\Gamma_k$, $k=2,\ldots,n$,
independent of $q^1$, for which
\[
\fpd{f}{q^k} = \Gamma_k(q^2,\ldots,q^n)\fpd{f}{q^1},\
k=2,\ldots,n\,.
\]
This is precisely the condition (\ref{conn}) from the introduction
(cf. Theorem 2). The connection $\A$ then reads $\A = dq^1 +
\Gamma_kdq^k$, with summation over $k =2, \ldots, n$. Note that
$Df = 0$ implies that the connection has vanishing curvature (the
horizontal distribution is involutive because it is annihilated by
an exact 1-form). Assume now that both of the above conditions hold
and keep the value of the momentum $\mu =
\partial L/\partial \dot q^1 - f$ fixed. We solve this
relation for $\dot q^1$ by writing $\dot q^1=\psi(q^k,\dot q^k)$,
with $k=2,\ldots,n$. The Routhian then is the function
\[R^\mu(q^k,\dot q^k)=L-(\mu+f)(\dot q^1+\Gamma_k\dot q^k),\] where
all instances of $\dot q^1$ on the right hand side have been
replaced by the function $\psi$. It now remains to compute the
2-form $\B_\mu^f$ which is the projection of
$d[(\mu+f)(dq^1+\Gamma_kdq^k)]$. After some straightforward
computations in which the condition $df^h=0$ is used, we obtain
\[
\B_\mu^f =
\frac12(\mu+f)\left(\fpd{\Gamma_k}{q^s}-\fpd{\Gamma_s}{q^k}\right)dq^k\wedge
dq^s.
\]
The latter is identically zero since the connection has zero
curvature due to $Df=0$. This also follows from the following
\begin{align*}
\fpd{\Gamma_k}{q^s}&= \frac{1}{\partial f/\partial q^1}\fpd{^2
f}{q^k\partial q^s} - \frac{1}{(\partial f/\partial
q^1)^2}\fpd{f}{q^k}\fpd{^2f}{q^s\partial q^1}\\ &= \frac{1}{\partial
f/\partial q^1}\fpd{^2 f}{q^k\partial q^s} - \frac{1}{\partial
f/\partial q^1}\Gamma_k\Gamma_s\fpd{^2f}{q^1\partial
q^1}=\fpd{\Gamma_s}{q^k}.
\end{align*}
We conclude that the Routhian reduction for Lagrangian systems with
a single quasi-cyclic coordinate is the Lagrangian system on the
reduced space with Lagrangian the Routhian $L-\A_\mu^f$. This
%shows
concludes the proof of
Theorem~\ref{thm:cyclic}.

\subsection{Functional Routhian reduction}\label{sec:funcrouth}
Our motivation for studying Routh-reduction for quasi-invariant
Lagrangians was inspired from the reduction technique called {\em
functional Routhian reduction} used in~\cite{funcrouth}. We will
argue here that functional Routhian reduction can be seen as
Routhian reduction for a quasi-invariant Lagrangian. Consider a
Lagrangian $L$ of type kinetic minus potential energy define on a
configuration space (locally) $(q^1,\ldots,q^{n-1},q^n)$. The
coordinate $q^n$ was denoted in~\cite{funcrouth} by $\phi$ and the
coordinates $q^k$ for $k=1,\ldots,n-1$ by $\theta^k$. The Lagrangian
$L$ is of the form
\[
L=\frac12 \left( M_{ij}(\theta) \dot q^i\dot q^j\right) -
W(\theta,\dot\theta,\phi) -V(\theta,\phi),
\]
with $M_{ij}(\theta)$ mass-inertia functions depending only on
$\theta^k$ and $W = (\lambda(\phi)/M_{nn}(\theta))
M_{nk}(\theta)\dot \theta^k$ and $V=V_{fct}(\theta) -
\frac12\lambda(\phi)^2/M_{nn}(\theta)$.

It should be immediately clear that $\phi$ is not a cyclic
coordinate, nor a quasi-invariant cyclic coordinate. We will however
define a `momentum map' $J_L^\lambda$ associated to the would-be
cyclic coordinate $\phi$:
\[
J^\lambda_L(\theta,\dot\theta,\phi)=\partial_{\dot \phi}
L(\theta,\phi,\dot\theta,\dot\phi) - \lambda(\phi) =
M_{kn}(\theta)\dot \theta^k + M_{nn}(\theta)\dot \phi
-\lambda(\phi).
\]
Note that, since $\lambda$ only depends on $\phi$ we may use the standard
connection $\A =d\phi$ when working in a local coordinate system.
The Lagrangian $L$ transforms as a quasi-invariant Lagrangian when restricted
to the level set $J^\lambda_L=0$:
\[
\left.\fpd{L}{\phi} \right|_{J^\lambda_L=0}= \lambda'(\phi)\dot\phi.
\]
Strictly speaking this example is not described in the theory
outlined above. We hope however that it is clear to the reader that
is an even more general type of Routh-reduction for quasi-invariant
Lagrangians that is valid only on a specific level set of the
momentum map. The correspondence between both techniques is also
seen from the fact that in~\cite{funcrouth} the authors define the
{\em functional Routhian} $L_{fct}$ as the function
\[
L_{fct}(\theta,\dot\theta)=\big(L (q^i,\dot q^i)-
\lambda(\phi)\dot\phi\big)_{J^\lambda_L=0}.
\]
This is precisely the function $\Ro^{\mu}$, with $\mu=0$ in our
analysis of quasi-invariant Lagrangians. Note that all regularity
conditions are satisfied and especially the horizontal condition
$dF_g^h=0$ is satisfied since $\lambda$ is independent of $\theta$.

\subsection{Charged particle in a constant magnetic field}
In~\cite{marmo88} the example of a charged particle in a constant
magnetic field $B$ is studied. The Lagrangian for this system is $L=
\frac12m(\dot x^2+\dot y^2)+eB(\dot xy-\dot yx)$. The coordinates
$x$ and $y$ are quasi-cyclic, and from this we may write that
$f(x,y)   = (-eBy,eBx)\in \R^2\cong\lag^*$. The 2-cocycle $\Sigma_f$
is not vanishing and proportional to $eB$. The (infinitesimal)
affine action on $\lag^*$ is completely determined by this 2-cocycle
$\Sigma_f$ and due to the abelian nature of the group, the
Lie-algebra of isotropy subgroup $G_{(\mu_1,\mu_2)}$ is trivial
since it is spanned by the kernel of $\Sigma_f$. In turn
$G_{(\mu_1,\mu_2)}=\{e\}$. The conserved momenta read: $m\dot
x+2eBy=\mu_1$ and $m\dot y -2eBx=\mu_2$. Therefore the quotient
space is $\R^2$. From the structure of the momenta equations it is
immediately seen that the system is $G$-regular. Further the
standard connection 1-form $\A = (dx,dy)^T$, with trivial horizontal
distribution implies that $Df=df^h=0$. Therefore all conditions are
met, and the Routhian is then a function on $\R^2$ depending only on
$x,y$:
\[
R^\mu=\frac{-1}{2m}\left((\mu_1-2eBy)^2 +(\mu_2+2eBx)^2\right).
\]
The symplectic 2-form $\B^f_\mu$ is precisely $2eBdx\wedge dy$. The
Routhian reduced equations of motion the read: $i_{(\dot x,\dot y)}
\B^f_\mu = dR^\mu$, or simply the momenta equations $m\dot
x+2eBy=\mu_1$ and $m\dot y -2eBx=\mu_2$.

{\bf ACKNOWLEDGMENTS}

J. V. is a Postdoctoral Fellow from the Research Foundation --
Flanders (FWO-Vlaanderen), and a Fulbright Research Scholar at the
California Institute of Technology. Additional financial support
from the Fonds Professor Wuytack is gratefully acknowledged. B.~L. is a part-time honorary postdoctoral researcher at the Department of Mathematical Physics and Astronomy, Ghent University, Belgium.

\end{document}